# La Emancipación Conceptual de Número Real de la Idea de Magnitud: Una Mirada Germánica

# The Conceptual Emancipation of Real Number from the Idea of Magnitude: A German Perspective.

MSc. Luis Giraldo González Ricardo
luis.gonzalez@matcom.uh.cu
Dr. Carlos Sánchez Fernández
csanchez@matcom.uh.cu
Universidad de La Habana

**Resumen**

En el presente artículo analizamos los aportes decisivos de tres miembros de la escuela alemana de matemática a la separación del concepto de número dela idea física de magnitud, durante la segunda mitad del siglo XIX. Además analizamos la importancia de la aritmetización del análisis para la consecuente axiomatización de los números reales, proceso que estos tres matemáticos promovieron en las obras principales que reseñamos.

**Abstract**

In the present article we study the decisive contributions of three members of the German mathematical school to the separation of the concept of number from the physic concept of magnitude, during the second half of the XIX[th] century. Besides we analyze the importance of the arithmetization of analysis to later axiomatization of real numbers, process that those three mathematicians promoted in the principal works we review.

## 1. Introducción.

Desde la Antigüedad Clásica las nociones de número y magnitud estuvieron íntimamente ligadas. Muestra de ello es la oposición por mucho tiempo de dar cabida dentro de la matemática a los números negativos y complejos, pues no satisfacían todas las reglas de la operatoria con magnitudes. Hasta bien entrado el siglo XIX los conceptos de magnitud y de número se perdían uno en el otro; en general el número era visto como una magnitud, de esta forma la fundamentación del concepto número se encontraba, por ejemplo, en la noción del tiempo.

En el presente trabajo pretendemos mostrar cómo en la obra de varios matemáticos de la escuela alemana se promueve la separación definitiva de los conceptos de número y de magnitud. Los matemáticos en cuestión son Leopold Kronecker (1823-1891), Richard Dedekind (1831-1916) y Georg Cantor (1845-1918).

Durante la primera mitad del siglo XIX se desarrolló un importante movimiento en búsqueda de un mayor rigor en el análisis, y en la matemática en general. Claros representantes de este movimiento son A. L. Cauchy, B. Bolzano y N. H. Abel, por solo citar tres.

Los conceptos básicos del análisis fueron puestos en tela de juicio, así aparecieron definiciones rigurosas de función continua, derivada en un punto y límite; este último se convirtió en el más básico de los tres y sirvió como fundamento del cálculo infinitesimal.

Poco a poco emergieron nuevos problemas y con ellos nuevos conceptos –como los de continuidad y convergencia uniforme- no se hicieron esperar. Pero los problemas se originaron no solo en las ramas superiores del análisis, también en sus raíces se encontraron cuestiones significativas que merecían revisión. Entre estas interrogantes esenciales, procurando una presentación rigurosa del concepto límite, apareció la necesidad de precisar el concepto número real.

Es a partir de la segunda mitad del siglo XIX que comienzan a aparecer distintas aproximaciones a la conceptualización de los números reales. De esta forma es que los personajes de nuestra historia aparecen en escena.

**2. Dedekind: el arquitecto de los números.**

Richard Dedekind nació en un pequeño pueblo del centro de Alemania: Braunschweig, el mismo lugar que vio nacer en 1777 a K. F. Gauss. Dedekind siguió los pasos de su predecesor y estudió en la Universidad de Göttingen, y allí defendió su doctorado en 1852 bajo la supervisión del propio Gauss.Dos años después de obtener su doctorado Dedekind presentó su *tesis de habilitación*, para poder enseñar en la universidad. El título de la memoria es **Über die Einführung neuer Funktionen in der Mathematik[1]**.Bien temprano se nota el particular interés de Dedekind por las operaciones, por encima de los elementos mismos. Preocupación que lo llevará a la visión abstracta de la matemática que logró en su madurez.

El objetivo principal de Dedekind en esta pequeña tesis es dar un método general para deducir la forma en que se deben extender las operaciones aritméticas de un dominio numérico más restringido a uno más general. Así lo describe el propio Dedekind:"*Esta conferencia no es sobre la introducción de una clase determinada de funciones dentro de la matemática (aunque se pudiera interpretar el título de esa forma); es más bien sobre la forma en que, en el desarrollo progresivo de la ciencia, nuevas funciones, o como se puede decir, nuevas operaciones son añadidas a las anteriores.*"[10, p. 755]

Desde el punto de vista matemático e histórico lo más interesante es el análisis que hace Dedekind del origen de las diferentes operaciones aritméticas. Considera que la más simple de ellas es la de encontrar el sucesor de cierto número; y esta repetida varias veces da lugar a la suma. Luego la repetición de la suma origina la multiplicación, y esta última da origen a la exponenciación.Otro detalle importante sobre las concepciones matemáticas de Dedekind que ya se destacan en esta temprana obra, es su noción de número, como ente creado por la mente humana. Dedekind plantea que los números naturales permiten la realización de las

---

[1]Sobre la introducción de nuevas funciones en la matemática.

operaciones de suma y de multiplicación sin problema alguno; no sucede lo mismo cuando se quieren realizar las operaciones indirectas [*indirekte*] o inversas [*umgekehrte*]: la resta y la división. La imposibilidad de realizar estas operaciones lleva a la *creación* de nuevos números donde sean válidas las operaciones indirectas: "*Y este requerimiento hace necesario crear una nueva clase de números, dado que con la sucesión original de los enteros positivos el requerimiento no se puede satisfacer.*"[10, p. 757]

En la tesis de 1854 Dedekind considera el origen de los números irracionales la operación de potenciación con exponente racional, idea que es desechada en su madurez, luego de que le preste más atención a la fundamentación de la aritmética de los números reales. Esta disertación de Dedekind contiene varias de las ideas desarrolladas en su madurez sobre la fundamentación de la matemática: la importancia de las operaciones realizadas sobre los elementos de cierto conjunto, por encima de los elementos en sí; además muestra de forma temprana su concepción más profunda sobre la filosofía de la matemática: la existencia de los entes matemáticos gracias a su creación por parte de la mente humana, en especial los números.

**2.1 Los números irracionales**

La principal motivación de Dedekind para emprender una construcción de los números irracionales, la encontró Dedekind en la docencia. En 1858 tuvo que enseñar los rudimentos del Cálculo Diferencial. Preparando dicho curso se enfrentó a que la fundamentación de varios teoremas básicos del Cálculo se basaban en argumentos geométricos, los que no consideraba válidos: era pertinente lograr una demostración puramente aritmética.

En **Stetigkeit und irrationale Zahlen**[2] Dedekind parte del conocimiento del dominio numérico racional **Q**, y muestra que existe una "operación" que no puede realizarse en este dominio: la operación de "cortar". Nótese como el problema de la continuidad de **Q** no se deriva de la imposibilidad de realizar una operación inversa, sino de una nueva. Pero la idea es semejante a lo establecido en su habilitación: se tiene una operación que no se puede realizar y es necesario aumentar el dominio numérico, de tal forma que sea posible dicha operación.

Véase cómo Dedekind comienza su memoria de 1872: "*Considero que la aritmética es una necesaria, o al menos natural, consecuencia del más simple acto, el de contar, y contar no es nada más que la creación sucesiva de la sucesión infinita de los enteros positivos, en el cual cada individuo es definido por el precedente, el simple acto de pasar de un elemento formado a uno por formar. (...) Mientras que la realización de estas dos operaciones[3] es siempre posible, la realización de las operaciones inversas ha sido en cada caso el motivo real de un nuevo acto creativo; entonces los números negativos y fraccionarios son creados por la mente humana; y con el conjunto de los números racionales se ha obtenido un instrumento de gran perfección.*"[11, p. 4]

Con esta cita se descubre cómo Dedekind mantiene buena parte de sus ideas presentes en 1854: del origen de los números naturales hasta la construcción de **Q** permanece inalterable. Pero desaparece la idea de la utilización de la potenciación con exponente no entero como origen de los números irracionales. Nótese también que en la fecha de publicación de esta monografía ya Dedekind había dedicado varios años al estudio del problema de la factorización única en los enteros algebraicos de una extensión simple de **Q**, que le llevó a la definición de cuerpo; lo que le permitió valorar mejor la importancia de estas operaciones básicas, que son características de un cuerpo[4].

---

[2] Continuidad y números irracionales.
[3] Adición y multiplicación. NA.

Se vio antes en **Über die Einführung…** cómo Dedekind ve en la necesidad de realizar nuevas operaciones la ampliación de los dominios numéricos. En **Stetigkeit…** no es ninguna de las operaciones comunes de la aritmética o el análisis, las que engendran los números reales, sino una nueva, a la que llama cortadura. Dedekind demuestra que la operación de corte es compatible con los números racionales, sin embargo, **Q** no es cerrado para esta operación. En otras palabras, existen cortes no definidos por números racionales. La idea de cortadura le proporciona los medios para caracterizar a continuidad de la línea recta, a la que le da validez de axioma:

"*Si todos los puntos de la línea recta se dividen en dos clases, tal que todo elemento de la primera clase se encuentra a la izquierda de todo elemento de la segunda clase, entonces existe un único punto que produce esta división de todos los puntos en dos clases, este corte de la recta en dos secciones.*"[11, p. 11]

En la memoria Dedekind se sirve de las cortaduras para la construcción de los números reales, y demostrar que cumplen el axioma de continuidad. Se debe hacer especial énfasis que Dedekind no menciona en ningún momento que esté desarrollando una teoría para las magnitudes irracionales, sino números y tampoco aparece ninguna mención a la noción de tiempo.

**2.2 Los números naturales.**

La memoria publicada por Dedekind en 1888 **Was sind und was sollen die Zahen?**[5] trata sobre la construcción de los números naturales, para lograr su objetivo evita la utilización de nociones geométricas o físicas, la metodología dedekindiana es, en esencia, teórico-conjuntista.

---

[4] Para detalles de los trabajos de Dedekind sobre los números algebraicos se puede consultar [15].
[5] ¿Qué son y para qué sirven los números?

En el mismo comienzo Dedekind deja clara la separación del concepto de número de las nociones de espacio y tiempo: "*Considero el concepto de número completamente independiente de las nociones o intuiciones de espacio y tiempo, porque lo considero un resultado inmediato de las leyes del pensamiento. Mi respuesta a los problemas propuestos en el título de esta monografía es, brevemente: los números son una creación libre de la mente humana; ellos sirven como medio para aprehender más fácilmente las diferencias entre las cosas. Es solo a través del proceso lógico de construir la ciencia de los números, y de esa forma adquirir el dominio numérico continuo, en que estamos preparados adecuadamente para investigar nuestras nociones de espacio y tiempo, poniéndolos en relación con este dominio numérico creado en nuestra mente*"[12, p. 31]

Nótese cómo Dedekind libera el concepto de número de la noción de tiempo, o en sentido más general, de la de magnitud. Antes de las construcciones de la segunda mitad del siglo XIX los números reales eran vistos como magnitudes, ejemplo fehaciente de ello se encuentra en la obra de Cauchy. Es decir, se utilizaron nociones físicas para explicar la continuidad de la recta real. Dedekind toma el camino inverso, los números reales sirven de modelo para las magnitudes físicas, dicho en otras palabras, es un medio para facilitar la comprensión del mundo físico. Además mantiene su principio de 1854, los números constituyen una creación de la inteligencia humana, aunque aclaró en [11] que la aritmética era consecuencia natural del acto de contar.

La construcción de Dedekind de **N**, está basada en ideas propias de la teoría de conjuntos: la existencia de un conjunto infinito destacado, con ciertas propiedades que pudiésemos llamar estructurales. Dedekind muestra en esta obra una mezcla armoniosa entre las ideas propias de su madurez, -utiliza aplicaciones e ideas conjuntistas-; con las líneas fijadas 30 años antes, definir la suma como la realización repetida de la operación sucesor.

### 3. La mirada analítica de Cantor.

La carrera matemática de Cantor comenzó siendo alumno de Kronecker, y a su sombra obtuvo su doctorado en la Universidad de Berlín. Una vez que llegó a Halle, bajo los auspicios de Eduard Heine (1821-1881), un discípulo de Dirichlet, se interesó por las rebeldes series trigonométricas. El estudio de la representabilidad de funciones por series trigonométricas era en ese momento una de las *vedettes* de la investigación matemática. La profundidad y belleza del tema conquistaron a Cantor. Las investigaciones por él desarrolladas eran sobre la unicidad de la representación de una función por serie trigonométrica, que era la extensión natural de un resultado previo de Heine. Cantor redujo el problema al análisis de la representatividad de la función idénticamente nula por una serie no nula divergente en cierto conjunto de puntos. En el desarrollo de su investigación Cantor necesitó de la caracterización de los números reales. Su construcción apareció en el artículo **Über die Ausdehnung eine Satzes aus der Theorie der trigonometrischen Reihen**[6] [5] de 1872.

Para su construcción de **R** Cantor utiliza las sucesiones de Cauchy, que tienen por característica esencial que no se necesita conocer el límite, para asegurar su "convergencia". Cantor utiliza esta propiedad para decir que cada sucesión fundamental tiene un límite definido $b$; y recalca que $b$ es solo un símbolo usado para destacar la naturaleza de la sucesión [5]. En el conjunto **B** de todos los símbolos $b$, Cantor demuestra que cada *elemento* es el límite de la sucesión que lo define. Luego de esto Cantor introduce un orden y las cuatro operaciones aritméticas. Para esto se basó en las propiedades de las sucesiones convergentes[7]. Una vez hecho esto Cantor comienza a llamar a los elementos de **B** como números. Según la interpretación de Dauben en [9] Cantor trata a los nuevos elementos de igual forma que a los números por conveniencia y simplificación, pues realmente son solo sucesiones de Cauchy.

---
[6]Sobre la extensión de un teorema de la teoría de las series trigonométricas.

Un detalle interesante de la construcción de Cantor es que aún se mantiene cerca de la tradición, pues utiliza para los nuevos números el término *Zahlgrösse* que significa *magnitud numérica*, en vez de *Zahl* que significa número.

Aunque se tiene definido un nuevo dominio numérico ¿será suficiente para completar a los racionales de forma continua? Para esto Cantor relaciona a los nuevos números con la recta. Una vez ubicado el origen, es sencillo relacionar las longitudes racionales con el número respectivo en **Q**; si la longitud puede aproximarse por una sucesión de longitudes racionales, tal sucesión es fundamental y define un número $b$ en **B**: "*La distancia a determinar desde el punto **b** al punto **O** es igual a **b**, donde **b** es el número determinado por la sucesión $\{a_n\}$.*" [9, p. 40]. En el sentido contrario no es tan sencillo, pues Cantor no tiene forma de asegurar la existencia de un punto en la línea para cada número; por ello utiliza la siguiente afirmación por axioma:"*(…) a cada número le corresponde un punto en la línea, cuya ordenada es igual al número.*" [9, p. 40]

Aunque Cantor identifica los números reales recién creados con los puntos de la recta, expresa que la existencia de los primeros es independiente de la recta; lo que dicho de otra forma, la existencia del número real es independiente de la geometría. Cantor prosigue en el camino de separar la geometría del análisis, además de independizar el concepto de número real del de magnitud. Para Cantor el número real es un punto, no una longitud; la que en esencia es una magnitud. A diferencia de Dedekind la separación realizada por Cantor es menos evidente, pero igualmente significativa.

---

[7] Para más detalles sobre la construcción de Cantor se puede ver [8] y [9].

En el año 1883 apareció una de las obras principales en el quehacer matemático de Georg Cantor: **Grundlagen einer allgemeinen Mannigfaltigkeitslehre[8]**. En ella aparecen los números transfinitos, un detallado estudio del continuo y además una buena parte de sus concepciones filosóficas. En este artículo Cantor considera a los números ordinales el origen de los números naturales, dados que los últimos surgen del acto de contar, o lo que es lo mismo, el acto de añadir unidades. Luego señala que los cardinales se obtienen como consecuencia del acto de contar, pues el número cardinal -*Anzahl*- de un conjunto es el número ordinal que se mantiene invariante tras cualquier reordenación.

En los **Grundlagen…** Cantor realiza un recuento sobre las teorías fundamentales de los números reales: la de Weierstrass, la de Dedekind y la propia. Cantor compara las tres construcciones poniendo sobre la mesa los pros y los contras de cada una. En el artículo se considera a Weierstrass como el primero en evitar el error lógico de definir los números irracionales por medio de límites de sucesiones con términos racionales Considera Cantor que la principal ventaja de las cortaduras es que, por encima de las otras dos es que a cada número real le corresponde una única cortadura, mientras que en las restantes un mismo número puede tener varias representaciones. Además G. Cantor realiza una exposición sucinta y más madura de su construcción de los números reales, en la que se refiere a los nuevos números recién creados como *Zahl* en vez de *Zahlgrösse*.

Finalmente en la sección §10 Cantor realiza un análisis detallado del continuo. El análisis comienza con una nota histórica sobre las diferentes concepciones del continuo, desde la Antigüedad Clásica, pasando por el Medioevo, hasta su época. Cantor considera que se ha

---
[8] Fundamentos de una teoría general de variedades.

tratado el continuo al nivel de dogma religioso, al ser considerado una intuición *a priori*, que incluso en su tiempo, es mal visto cualquier intento de estudiar el continuo desde la matemática, utilizando la razón: "*Cualquier intento aritmético de determinación de este misterio –el del continuo- es visto como una invasión prohibida y es rechazada con vigor. Mentes tímidas entonces, se llevan la impresión de que el continuo no es un problema de un concepto matemático lógico, sino más bien un dogma religioso.*"[6, p 903]

A continuación Cantor critica las concepciones anteriores de considerar el tiempo y el espacio modelos de continuidad, lo que considera no admisible. Lo anterior se debe a que la continuidad del tiempo está condicionada por la consideración del continuo asumido de forma independiente al tiempo; luego no puede ser interpretado desde el punto de vista objetivo ni subjetivo: "*(…) en mi opinión utilizar la intuición del tiempo para discutir el concepto más general de continuo no es la forma correcta de proceder; el tiempo es en mi opinión una representación.*"[6, p. 904]. Además: "*Es también mi convicción de que con la forma de intuición del espacio, no se puede ni tan siquiera comenzar a adquirir conocimiento del continuo*" [6, p. 904]

En particular Cantor realiza su análisis en el espacio *n*-dimensional $\mathbf{R^n}$ con la distancia euclidiana. Apoyado en las ideas topológicas para el estudio de los conjuntos de puntos de $\mathbf{R}$ –punto de acumulación, conjunto derivado- encontró que una condición necesaria para la continuidad de un conjunto de puntos era que este tiene que coincidir con su conjunto derivado, lo que denominó *perfecto*[9]. Rápidamente Cantor da un ejemplo de un conjunto perfecto que no es continuo: el conjunto conocido hoy en día como *conjunto de Cantor,* que es no numerable, perfecto, no contiene puntos interiores y claramente no es continuo. Para

---

[9]En terminología actual sería cerrado.

completar la caracterización del continuo Cantor considera que falta el concepto de *conexidad*:

"*Llamamos a* T *un conjunto conexo si, para cualesquiera dos puntos* t *y* t'*, y para cualquier número arbitrariamente pequeño α existe siempre una cantidad finita de puntos* $t_1, t_2,...,t_v$ *tal que las distancias* $tt_1, t_1t_2,...t_vt'$ *son todas menores que α*" [6, p. 906]

A la luz de los métodos actuales la conexidad cantoriana es bastante extraña, pues siguiendo a J. Ferreirós en [18], esta definición de conexidad contempla a **Q** como un conjunto conexo; realmente la idea de Cantor está más ligada a la densidad. Sin embargo, muestra la utilidad de la topología en el camino a la caracterización del continuo.

De esta forma Cantor caracterizó a los conjuntos continuos como "*conjuntos perfectos y conexos*" [6, p. 906], y en particular **R** satisface esta definición; de la que se entiende que los números reales para Cantor son puntos, no longitudes. Finalmente Cantor concluye realizando una crítica a los enfoques de Bolzano y Dedekind. Del primero señala que en **Paradoxen der Unendlichen** [10] solo resalta la característica que tienen los conjuntos continuos de ser conexos; mientras que del segundo dice que en **Stetigkeit…** solo destacó la característica común a todos los conjuntos continuos de ser cerrados.

Sobre la crítica a Dedekind se puede decir que el objetivo de su memoria de 1872, no era caracterizar la continuidad del espacio en general, sino de los números reales en particular. Mientras que Cantor, dado su interés particular en los problemas de numerabilidad, teoría de conjuntos de puntos y topología, se centra en buscar una caracterización más apropiada a sus necesidades; por ello la caracterización de Cantor utiliza las propiedades métricas del espacio *n*-dimensional. Por otro lado Dedekind utiliza la estructura de orden de los números racionales.

---

[10]Paradojas del infinito.

Aunque con notables diferencias, las investigaciones de Cantor y Dedekind comparten el mismo principio básico: el continuo aritmético es cierto conjunto de **puntos** que satisfacen cierta propiedad, en cada caso la caracterización de la continuidad. Ahora, ¿es esta la única forma de considerar la continuidad de los números reales?

**4. Kronecker, una visión discordante.**

Los escritos de Kronecker sobre los fundamentos de la matemática no son abundantes, además de existir numerosas leyendas sobre él dentro de la comunidad matemática[11]. En el referido artículo aparece la cita de una carta a Cantor de 1884, donde de modo general explica parte de sus ideas sobre la forma de hacer matemática:

"*Porque tú tomaste mis cursos hace más de 20 años, y has tenido contacto casi ininterrumpidamente conmigo desde entonces, has escuchado mis posiciones lo suficiente como para comprender mejor de lo que pudiera yo describir aquí, y que -habiendo profundizado muy temprano en cuestiones filosóficas bajo la tutela de Kummer- reconocí, al igual que Kummer hizo, la inestabilidad de tales especulaciones y tomé refugio en el paraíso de la matemática real. (...) en mi obra matemática he tenido mucho cuidado de expresar sus fenómenos y verdades en una forma lo más libre posible de conceptos filosóficos, y consecuentemente he comenzado el camino para basar toda la matemática en la teoría de los números naturales, y creo que esto puede hacerse sin excepción. Admito que este punto es solo mi creencia. Pero donde ha tenido éxito, veo el verdadero progreso, incluso -o porque- es una regresión a los principios más simples, más aún porque prueba que los nuevos conceptos introducidos son al menos innecesarios. (...) Reconozco un verdadero valor*

---

[11]En este punto se recomienda la lectura de [16].

*científico -en el campo de la matemática- solo en verdades concretas, o para ponerlo más claro, solo en fórmulas matemáticas."* [16, p. 45]

Entre líneas se puede leer la forma en que Kronecker pretende utilizar el concepto de número natural como fundamento de la matemática. La idea básica es aritmetizar la matemática -pura-, pero no sobre la base de los números reales, sino en el idioma más básico de los números naturales.

Kronecker publicó sus ideas en el artículo **Über den Zahlbegriff**[12] [21] de 1887. Este artículo y su ampliación posterior –también de 1887-, junto al curso del semestre de verano de 1891([22], [1]) son las únicas fuentes de la obra publicada de Kronecker, útiles para obtener información sobre sus posturas referentes a la aritmetización y el concepto número.

En la introducción del **Zahlbegriff** Kronecker realiza una separación entre las principales disciplinas de la matemática, que según él son: la geometría, la mecánica y la aritmética. Conjuntamente a esta división de la matemática Kronecker expone su objetivo: *"Creo que nosotros tendremos éxito en 'aritmetizar' todo el contenido de estas disciplinas[13], esto es, fundamentarlas solamente en el concepto número tomado en su sentido más estrecho, y entonces desechar las modificaciones y extensiones de este concepto[14], las cuales son ocasionadas por las aplicaciones a la geometría y a la mecánica."* [21, p. 949]

Nótese cómo Kronecker considera a las cantidades irracionales fuera de los dominios de la aritmética. Kronecker estima que la noción de número irracional es propia de la geometría y la mecánica, no de la aritmética. En la separación que realiza Kronecker de las disciplinas matemáticas fundamentales, la aritmética no es solo la rama de la matemática que estudia las propiedades de los números, según Kronecker: *"La palabra aritmética no debe*

---

[12]Sobre el concepto de número.
[13]Áritmética, álgebra y análisis. NA
[14]Me refiero específicamente a la adición de las cantidades irracionales y las cantidades continuas. (Nota de Kronecker.)

*comprenderse en el sentido restringido que es usual, sino incluyendo todas las disciplinas matemáticas, -con excepción de la geometría y la mecánica- especialmente el álgebra y el análisis."* [21, p. 949]. Esta es la diferencia principal con los otros aritmetizadores alemanes: Dedekind, Cantor, Heine y Weierstrass; quienes consideraban a los números irracionales como elementos propios e imprescindibles dentro del análisis.

Anteriormente se vio cómo Dedekind define los números naturales, a diferencia de Kronecker que comienza su estudio a partir de ellos, más precisamente los números ordinales: *"Los números ordinales constituyen de forma natural el punto de partida para el desarrollo del concepto número"* [21, p. 949]. Es en este contexto que surge la pregunta siguiente: ¿Por qué los números ordinales como origen de sus estudios? Para esto hay que adentrarse un poco más en las concepciones de Kronecker sobre los fundamentos de la matemática. Donde más detalles existen sobre las ideas de Kronecker es en los manuscritos del curso de 1891, que lleva, precisamente, el título **Über den Begriff der Zahlen in der Mathematik**[15][1].

En el curso Kronecker establece cuatro principios fundamentales sobre la forma de concebir la matemática [1, pp. 210,232-233]:1) La disciplina matemática no admite sistematización; 2) La matemática debe ser tratada igual que las ciencias de la naturaleza, puesto que sus objetos son tan reales como los de sus ciencias hermanas; 3) Es necesario evitar la mezcla de una disciplina con la otra; 4) Utilizar en cada disciplina un método acorde a la materia que se estudia.Lo planteado por Kronecker no quiere decir que la matemática sea una ciencia natural, sino que debe ser tratada como tal. Para Kronecker la matemática es una ciencia pura, pero debe, al igual que las ciencias naturales, obtener las verdades inobjetables de la experiencia. El único elemento que la experiencia ha dado a la matemática es el acto de;

---

[15]Sobre el concepto de número en la Matemática.

contar, por ello los números ordinales son la única verdad incuestionable de la matemática. Luego de este análisis cobra total sentido la célebre cita que H. Weber atribuye a Kronecker: "*Dios hizo los números naturales y el resto es obra de los hombres*".

Con los números ordinales por punto de partida, Kronecker define los números cardinales, y estos últimos son los que se entienden por números. Grosso modo, Kronecker utiliza la biyección entre todos los conjuntos con la misma cantidad de elementos. El método es el siguiente, si se tiene una colección de objetos el número cardinal que le corresponde, es el último número ordinal utilizado para contar la cantidad de elementos -*Anzahl*-. La idea de Kronecker se asemeja a la utilizada por Russell y Frege para definir el concepto de número; es necesario señalar que la principal diferencia es que para Kronecker los números ordinales **ya están dados**, mientras que Frege y Russell los están definiendo a partir de conceptos lógicos. Es necesario aclarar que Kronecker no define al número cardinal $n$ como la clase de todos los conjuntos con $n$ elementos, sino que utiliza a un representante de la clase, que en este caso es el número ordinal $n^{mo}$. Además, -utilizando lenguaje matemático actual- Kronecker demuestra que la operación de asignarle a cada conjunto su correspondiente número cardinal está bien definida, es decir, que es independiente del orden de los elementos. La definición de las operaciones directas básicas, la suma y la multiplicación, es en esencia la misma que escogió Dedekind. En efecto, la suma es la repetición del acto de contar, y la multiplicación la repetición del acto de sumar.

La parte más interesante del trabajo de Kronecker está en el más corto de sus epígrafes: **Cálculo con Variables**. Su idea se basa en la introducción de indeterminadas, lo que le permitirá eliminar la noción de número irracional algebraico, incluso los números negativos.

Por ejemplo, dice Kronecker que la igualdad 7-9=3-5 se puede sustituir por la congruencia módulo *(x+1)*: *7+9x≡3+5x* mod *(x+1)*.

En este corto artículo Kronecker no cumple todos sus objetivos, pues no expone la forma en que se puede omitir la introducción de las irracionalidades algebraicas, y solo tratarlas a partir de su respectivo polinomio irreducible. Tal labor la completó en la ampliación publicada en el *Journal de Crelle*. La ampliación se reduce a la sección **Cálculo con variables**, donde Kronecker muestra la forma de extender a los números racionales y a los irracionales algebraicos el método utilizado con los números negativos.

La primera adición que realiza Kronecker es la introducción de una forma de tratar a los números racionales como indeterminadas [22, p. 6]. Más adelante recalca que "*De este modo la sucesión de fracciones racionales no se define solamente, sino queda justificada también.*" [22, p. 6]. De esta forma Kronecker incorpora a los números racionales a los dominios de la aritmética, al interpretarlos como congruencias con indeterminadas sobre el conjunto de los números naturales. De esta forma justifica el tratamiento de los racionales dentro de la aritmética propiamente dicha. También se puede notar que no es su único objetivo fundamentar la matemática en el concepto de número, sino trasladar los métodos propios de la aritmética, congruencias, a otras ramas.

La otra extensión realizada por Kronecker es el tratamiento de las irracionalidades algebraicas, de forma similar a lo hecho con los números negativos y los racionales. La extensión provista por Kronecker es más trabajosa que en los casos anteriores, pues las expresiones que sirven para definir las irracionalidades algebraicas son de grado mayor al 2; además de que una ecuación irreducible no solo define un número algebraico, sino a sus conjugados también. El método seguido por Kronecker para evitar la utilización de las irracionalidades algebraicas es mucho más complicado que para los casos anteriores. En sus trabajos donde extiende los resultados de Kummer sobre la factorización única, Kronecker había evitado la introducción de los números algebraicos, pero para ello solo necesitó el hecho de que un número algebraico α sea raíz de cierto polinomio irreducible *p(x)* con coeficientes en **Q**. En este caso necesita el valor de α, por eso es más complicada la solución. En lo siguiente veremos la idea general de lo realizado por Kronecker. Consideremos *p(x)* un polinomio irreducible sobre **Q**, el objetivo de Kronecker es encontrar un número *n* en **N**, tal que le permita aislar las raíces reales de *p(x)* en intervalos de longitud *1/n*. Se entiende por aislar las raíces que en cada intervalo de longitud *1/n* tal que en los extremos *p(x)* tome valores con signos distintos, haya solo una raíz; y de lo contrario, no tenga raíz alguna. Por otro lado, los cálculos con el número α se reducirán a cálculos con aproximaciones racionales, cada vez más precisas.

Compruébese con las palabras de Kronecker: "*La llamada existencia de las raíces irracionales de las ecuaciones algebraicas está fundamentada únicamente en la existencia de intervalos con la característica especificada; la legitimidad de calcular con las raíces individuales de una ecuación algebraica está basada completamente sobre la posibilidad de aislarlas, por tanto, en la posibilidad de determinar un número, el cual designamos por* s *arriba[16]. Si tal número* s *está determinado por la propiedad de que los intervalos de longitud* 1/s *son suficientemente pequeños como para aislar las distintas raíces de la misma ecuación, entonces **mayor que** y **menor que** se definirán simplemente a través de la sucesión de los intervalos de aislamiento. (…) La verdadera esencia del asunto, sin embargo, queda completamente clara cuando en la deducción anterior se evita la utilización de las fracciones y se hace uso exclusivo de los números naturales.*" [22]

---

[16]En la memoria de Kronecker.

El procedimiento empleado por Kronecker es bien característico de su forma de hacer matemática. Kronecker logra encontrar el *s* en **N** con las condiciones requeridas de forma constructiva y solamente utilizando los coeficientes del polinomio *p(x)*. Es decir, Kronecker da un procedimiento para lograr aproximaciones cada vez más precisas para los irracionales algebraicos, por medio de los coeficientes de *p(x)*, que pueden incluso considerarse enteros.

Ni en la memoria original, ni en la versión aumentada Kronecker hace referencia a algún método que proporcionara un resultado similar para los irracionales trascendentes, lo que está en consonancia con lo expuesto al inicio del artículo, donde considera al continuo ajeno a los dominios de la aritmética. Realizando una interpretación exegética de este trabajo se puede decir que para Kronecker la aritmetización consistía en fundamentar la matemática en los números naturales, sin llegar a la obtención del continuo, del modo en que lo hicieron Dedekind y Cantor.

Finalmente, en el ya referido curso de 1891 Kronecker no hace referencia alguna al tratamiento de los números algebraicos. Pero poniendo la situación en contexto, es bastante lógico, puesto que en 1888 publicó un extenso trabajo dedicado a exponer sus generalizaciones de los resultados obtenidos por Kummer. El título de la obra es **Zur Theorie der allgemeinen complexen Zahlen und der Modulsysteme**[17].

Por supuesto que en este momento vale preguntarse por qué los números irracionales no son parte de la aritmética. La respuesta debe, nuevamente, buscarse en las tesis sobre los fundamentos esgrimidas por Kronecker en su curso de 1891.La tercera tesis refiere que es necesario, o mejor, imprescindible, evitar la mezcla entre las diferentes disciplinas de la matemática; y la cuarta refiere que cada disciplina debe desarrollar sus propias herramientas.

---

[17]Sobre la teoría de los números complejos generales y de los sistemas de módulos.

Dicho de otra forma, lo que es inherente a la geometría debe permanecer en ella, y no mezclarse con la aritmética.

El porqué Kronecker consideraba los números irracionales pertenecientes solo al campo de la geometría parece no estar claro, pero se puede establecer una conjetura bastante plausible, y está relacionada con el origen de los números irracionales. Los números irracionales, o mejor dicho, las magnitudes inconmensurables se dieron a conocer por medio de la geometría helénica. En nuestra opinión esta es la razón que consideramos más acorde para explicar las ideas de Kronecker. Es consecuencia inmediata de sus tesis que el número irracional no podía ser parte de la aritmética. Luego surge una interrogante bastante obvia, y que Kronecker dejó sin responder: si el análisis es parte de la aritmética, entonces tampoco se puede aceptar la utilización de los números irracionales en la fundamentación del análisis, entonces, ¿cómo fundamentar el concepto de límite sin los números irracionales? En el curso de 1891 solo arroja un poco de luz sobre las aplicaciones del análisis, en este caso solo bastaría conocer aproximaciones racionales de cierta cantidad irracional [1].

Igualmente es importante destacar que Kronecker habla de cantidades irracionales, no de números irracionales. Este hecho se debe a que, según Kronecker, los irracionales no son números, pues este honor solo pertenece a los números naturales. Si se considera como concepto primario el de número natural, que tiene su origen en el acto de contar; entonces, para incluir a los enteros y a los racionales dentro de los números, debemos expandir el concepto primario de número. La ampliación del concepto se debe a que ni los ``números" enteros ni los racionales sirven para contar. Claramente el objetivo de Kronecker, al desarrollar herramientas para evitar el uso de los números racionales, por ejemplo, es evitar la generalización del concepto de número, pues este debe permanecer inalterable. Esta afirmación también es compartida por J. Boniface en [3].

Por supuesto que estas ideas lo ponen en desacuerdo con otros aritmetizadores, y el caso más claro es Dedekind. Se vio anteriormente que Dedekind desde 1854, se preocupa por la expansión de los conceptos ya existentes, en particular la ampliación del concepto de número. No es esta la única contradicción de Kronecker con las otras posturas dentro de la aritmetización. Para Kronecker las construcciones de Cantor y Dedekind de los números reales carecían de sentido, puesto que ambas utilizaban conjuntos infinitos de números racionales de una forma u otra. La utilización por parte de Dedekind y Cantor del infinito actual para sus construcciones de los números reales no concordaban con las ideas constructivistas de Kronecker, que entre otras cosas, no consideraba válido el infinito actual.

Se ha podido constatar que el proyecto de Kronecker estaba dirigido a fundamentar la aritmética en bases más internas. Kronecker toma por base el concepto de número ordinal, y con este construye su aritmética, evitando las generalizaciones del concepto de número. Sin duda el proyecto de Kronecker contiene algunas lagunas, así es el caso de la fundamentación del concepto de límite sin números irracionales en el análisis.

5. **Conclusiones.**

La metodología y concepciones de Cantor, Kronecker y Dedekind sobre la actividad matemática son divergentes en la mayoría de los casos. Pero no solo hay divergencias en los enfoques de estos tres matemáticos, también existen puntos de contacto. Las diferencias son más notables entre Cantor y Dedekind por un lado, y Kronecker del otro. Los dos primeros adoptaron una posición teórico-conjuntista, en su versión inicial y más ingenua, mientras que Kronecker no aceptaba la utilización de conjuntos arbitrarios.

Los métodos finitistas de Kronecker no permitían la utilización del infinito real, como objeto matemático bien definido; en ese sentido Kronecker permaneció ligado a la idea del infinito potencial, más una forma de hablar que una herramienta matemática. Dedekind aceptó el infinito en su quehacer matemático, tanto los números reales, los módulos y los ideales estaban definidos por medio de conjuntos infinitos de elementos que comparten cierta propiedad. Cantor no solo utilizó conjuntos infinitos para su caracterización de los números reales, sino que dedicó gran parte de su vida a desentrañar los misterios de la infinitud. Sobre este apartado Kronecker escribió: "*Las consideraciones anteriores me parecen opuestas a la introducción de aquellos conceptos de Dedekind, digamos 'módulo', 'ideal', etc., de modo similar (ellas son opuestas a) la introducción de varios conceptos con cuya ayuda en tiempos recientes y con frecuencia se ha intentado (...) concebir y establecer los 'irracionales' en general.*" [9]

Cantor y Dedekind creían que la actividad matemática era puramente lógica, es decir, la matemática es una actividad intelectual; tal concepción se refleja en la libertad de creación que ambos matemáticos desarrollaron en su quehacer profesional. Además, ya antes se señaló cómo Dedekind consideraba a los números una creación de la mente humana; mientras que Cantor desarrolló toda una teoría sobre los números transfinitos, los que constituyeron una creación esotérica de su intelecto. Por el lado contrario Kronecker veía a la matemática como una ciencia basada en la experiencia, y por tanto que se descubre; luego la libertad de creación está limitada por lo que la experiencia pueda aportar. Aunque ambos matemáticos se interesaron por la caracterización de la continuidad, existen diferencias en sus enfoques y objetivos particulares: 1) Dedekind se interesó por la caracterización de la continuidad con el objetivo de alcanzar una formulación precisa del concepto número real; 2) El interés de

Cantor de identificar las características esenciales del continuum era con la finalidad de probar la hipótesis del continuo; 3) La mirada de Cantor al continuo fue típica de un analista, por lo que se centró en sus propiedades métricas y topológicas; 4) La aproximación seguida por Dedekind estuvo más influida por su tendencia "estructuralista", por lo que su interés se concentró en la relación de orden que se establece entre los números reales a partir de los racionales; y además cómo se transmiten las operaciones aritméticas[18].

Por otro lado, en la obra de Cantor se destaca la necesidad de estudiar el continuo desde el punto de vista topológico, no solo desde la aritmética. Como bien se señala en [Sin92], no es evidente reconocer el lenguaje topológico en las propiedades utilizadas para caracterizar la continuidad de **R**, dígase: 1) Propiedad de las cortaduras o axioma de Dedekind; 2) Convergencia de las sucesiones de Cauchy; 3) Teorema de Bolzano de los valores intermedios; 4) Principio de intervalos encajados; 5) Principio del supremo. Especialmente por las diferentes propiedades estructurales que posee el conjunto de los números reales: es un cuerpo ordenado desde el punto de vista algebraico, además de arquimediano y "saturado" o completo en el sentido topológico asociado a la noción usual de distancia. Es relevante destacar que Dedekind en 1872 había revelado las propiedades de orden de los números reales, puesto que las utiliza de forma fundamental en su construcción. El estudio del continuo por Dedekind es solo desde el punto de vista aritmético, mientras que Cantor resalta la necesidad del estudio topológico.

Un punto en común entre estos tres matemáticos es su visión aritmetizadora de la matemática, cada uno desde su perspectiva, pero los tres percibieron a la aritmética como fuente definitiva de rigor; por ello la fundamentación de la matemática debía buscarse en la aritmética. En ellos se establece un punto de inflexión en la caracterización del concepto de número, al separarlo de la naturaleza geométrica que había tenido hasta ese momento, con ellos -en gran medida- se deja de centrar el discurso sobre la noción de número en el concepto magnitud.

---

[18]Sobre las tendencias estructuralistas de Dedekind y su relación con el estudio del continuo aritmético se puede ver [19]

Para comprender mejor lo que se ha expuesto en el párrafo anterior, se puede analizar la posición de A. L. Cauchy en su celebrado **Cours d'Analyse**[19] de 1821. En los *Preliminares* Cauchy aclara varios conceptos básicos para el desarrollo de su texto, entre ellos aclara las diferencias entre **número** y **cantidad**: "*Indicaremos qué ideas serán más apropiadas para adjuntarle a las palabras 'número' y 'cantidad'. Tomaremos el significado de 'número' en el sentido que es usado en aritmética, donde los números surgen de la medida absoluta de magnitudes, y aplicaremos solamente el término cantidades reales positivas o negativas a los números precedidos por los signo + o -. Además reconoceremos que estas magnitudes tienen el propósito de expresar crecimiento o decrecimiento; entonces una magnitud será representada por un número solo si queremos compararla con otra magnitud del mismo tipo tomada como unidad (...)*" [4, p. 5]

Para Cauchy los números son solo los reales positivos, mientras que los reales negativos son cantidades, es decir, no son números. La posición de Cauchy es muestra de la confusión entre número y magnitud de inicios del siglo XIX. El concepto de número en Cauchy no es lo suficientemente amplio, de modo que incluyese a los números negativos, lo que puede considerarse un punto de contacto con Kronecker. Sin embargo, Cauchy considera que el número surge de medir, por tanto no considera a los números naturales el punto de partida para el estudio de la aritmética. Véase la diferencia con Cantor, Dedekind o Kronecker.

---

[19]Curso de Análisis

El hecho de independizar al concepto de número de nociones ajenas a su naturaleza aritmética, como lo es el tiempo o la longitud, constituye un momento determinante hacia la axiomatización del concepto de número. Lo anterior significa que la esencia del concepto comienza a mostrarse de forma pura, es decir, se desliga de nociones ajenas para obtener sus características esenciales, las que sirven de base a sus axiomas. Esta es una fase posterior digna de otro estudio semejante, por ejemplo, a través de la obra de figuras germánicas del calibre de David Hilbert y Hermann Weyl.

**Bibliografía.**